\newtheorem{theorem1}{Theorem}

\newcounter{theor}
\newtheorem{theorem}{Theorem}[theor]
\newtheorem{lemma}{Lemma}[theor]
\newtheorem{corollary}{Corollary}[theor]

\newenvironment{proof}[1][Proof]{\begin{trivlist}
\item[\hskip \labelsep {\bfseries #1}]}{\end{trivlist}}

\documentclass[12pt]{article}

\begin{document}
\centerline{\bf \large THE $\alpha$-INVARIANT ON $\bf
CP^2\#2\overline{CP^2}$}
\bigskip
\centerline{Jian Song}
\bigskip

\centerline{Department of Mathematics} \centerline{Columbia
University, New York, NY 10027}
\bigskip
\bigskip

\centerline{\bf \S 1. Introduction}
\bigskip

\noindent The global holomorphic invariant $\alpha_{G}(M)$
introduced by Tian \cite{T1}, Tian and Yau \cite{TY} is closely
related to the existence of K\"ahler-Einstein metrics. In his
solution of the Calabi conjecture, Yau \cite{Y} proved the
existence of a K\"ahler-Einstein metric on compact K\"ahler
manifolds with nonpositive first Chern class. K\"ahler-Einstein
metrics do not always exist in the case when the first Chern class
is positive, for there are known obstructions such as the Futaki
invariant. For a compact K\"ahler manifold $M$ with positive Chern
class, Tian \cite{T1} proved that $M$ admits a K\"ahler-Einstein
metric if $\alpha_{G}(M)>\frac{n}{n+1}$, where $n=\dim M$. In the
case of compact complex surfaces, he proved that any compact
complex surface with positive first Chern class admits a
K\"ahler-Einstein metric except $CP^2\#1\overline{CP^2}$ and
$CP^2\#2\overline{CP^2}$ \cite{T3}. It would be also interesting
to find the estimate of the $\alpha$ invariant for
$CP^2\#1\overline{CP^2}$ and $CP^2\#2\overline{CP^2}$. In this
paper, we apply the Tian-Yau-Zelditch expansion of the Bergman
kernel on polarized K\"ahler metrics to approximate
plurisubharmonic functions and compute the $\alpha$-invariant of
$CP^2\#2\overline{CP^2}$. This gives an improvement of
Abdesselem's result \cite{A}. More precisely, we shall show that:
\bigskip
\begin{theorem1}
{$\alpha_{G}(CP^2\#2\overline{ CP^2})$=$\frac{1}{3}$.}
\end{theorem1}

\bigskip

\noindent Let ($M$, $\omega$) be a compact K\"ahler manifold,
where $\omega$=$\sqrt{-1} g_{i \overline{j}}dz_i \wedge d
\overline{z}_j$. We will also prove Tian's conjecture on the
generalized Moser-Trudinger inequality in the special case where
$\alpha_{G}(M) > \frac{n}{n+1}$, for $n=\dim M$. Let
$$P(M,\omega)=\left\{\phi\;|\;\omega_{\phi}=\omega+\partial\overline{\partial}\omega>0,\;\sup_{M}\phi=0\right\}.$$
Let $F_{\omega}$ and $J_{\omega}$  be the functionals defined on
$P(M,\omega)$ by
$$ F_{\omega}(\phi)=J_{\omega}(\phi)-\frac{1}{V}\int_{M}\phi\omega^n-
   \log(\frac{1}{V}\int_{M}e^{h_{\omega}-\phi}\omega^n)$$
$$ J_{\omega}(M)=\frac{1}{V}\sum\limits^{n-1}_{i=0}\frac{i+1}{n+1}
   \int_{M}\partial\phi\wedge\overline{\partial}\phi\wedge\omega^i\wedge\omega_{\phi}^{n-i-1}.
$$

\noindent Assume $\left( M,\omega_{KE} \right) $ is a
K\"ahler-Einstein manifold with positive first Chern class and
$Ric(\omega_{KE})=\omega_{KE}$, then for any $\phi \in P(M,\omega
_{KE})$, Ding and Tian \cite{DT} proved the following inequality
of Moser-Trudinger tpye:
$$
\frac{1}{V}\int_{M}e^{-\phi }\omega ^{n}\leq Ce^{J_{\omega }(\phi )-\frac{1}{%
V}\int_{M}\phi \omega ^{n}}.$$

\noindent Tian\cite{T4} also conjectured that
$\frac{1}{V}\int_{M}e^{-\phi }\omega ^{n}\leq
Ce^{(1-\delta)J_{\omega }(\phi )-\frac{1}{ V}\int_{M}\phi \omega
^{n}}$ for $\delta > 0$ sufficiently small, if $\phi$ is
perpendicular to $\Lambda _{1}$, the space of eigenfunctions of
$\omega_{KE}$ with eigenvalue one.

\noindent We shall prove:

\begin{theorem1}
Let \ $\left( M,\omega \right) $ be a K\"ahler manifold with
positive first Chern class. Assume that $\alpha
(M)>\frac{n}{n+1},$ so that $M$\ admits a K\"ahler-Einstein metric
$\omega _{KE},$ and \ there exist \ constants $ \delta =\delta
\left( n,\alpha (M)\right) $\ \ and \ $C=C\left( n,\lambda
_{2}(\omega_{KE})-1,\alpha (M)\right)$ such that for any $\phi \in
P(M,\omega _{KE})$ which satisfies $\phi \perp \Lambda _{1}$, we
have:\
$$
\ F_{\omega _{KE}}(\phi )\geq \delta J_{\omega _{KE}}(\phi )-C
$$
Here $\lambda_{2}(\omega_{KE})$ is the least eigenvalue of
$\omega_{KE}$ which is bigger than 1.
\end{theorem1}

\bigskip
\bigskip
\noindent $Acknowledgements$. The author deeply thanks his
advisor, Professor D.H. Phong for his constant encouragement and
help. He also thanks Professor G. Tian for his suggestion on this
work. This paper is part of the author's future Ph.D. thesis in
Math Department of Columbia University.

\bigskip
\bigskip

\centerline{\bf \S 2. Holomorphic approximation of psh}
\bigskip
\setcounter{theor}{2} \noindent In this section, we will employ
the technique in \cite{T2,Z} to obtain the approximation of
plurisubharmonic functions by logarithms of holomorphic sections
of line bundles. The Tian-Yau-Zelditch asymptotic expansion of the
potential of the Bergman metric is given by the following theorem
\cite{T2,Z}.
\bigskip
\begin{theorem} Let M be a compact
complex manifold of dimension n and let $(L,h)\rightarrow M$ be a
positive Hermitian holomorphic line bundle. Let g be the K\"ahler
metric on M corresponding to the K\"ahler form $\omega
_{g}=Ric(h)$. For each $m\in N$, $h$ induces a Hermitian metric
$h_{m}$ on $L^{m}.$ Let
$\{S_{0}^{m},S_{1}^{m},...,S_{d_{m-1}}^{m}\}$ be any orthonormal
basis of $H^{0}(M,L^{m})$, $d_{m}=\dim H^{0}(M,L^{m}),$ with
respect to the inner product:
$$
(S_{1},S_{2})_{h_{m}}=\int_{M}h_{m}(S_{1}(x),S_{2}(x))dV_{g},
$$
where $dV_{g}=\frac{1}{n!}\omega _{g}^{n}$ is the volume form of
$g$. Then there is a complete asymptotic expansion:

$$
\sum\limits_{i=0}^{d_{m}-1}||S_{i}^{m}(x)||_{h_{m}}^{2}=a_{0}(x)m^{n}+a_{1}(x)m^{n-1}+a_{2}(x)m^{n-2}+...
$$
for some smooth coefficients $a_{j}(x)$\ with $a_{0}=1$. More
precisely, for any k:
$$
\left| \left|
\sum\limits_{i=0}^{d_{m}-1}||S_{i}^{m}(x)||_{h_{m}}^{2}-
\sum_{j<R}a_{j}(x)m^{n-j}\right| \right| _{C^{k}}\leq
C_{R,k}m^{n-R}
$$
where $C_{R,k}$\ depends on $R,k$\ and the manifold $M$.
\end{theorem}

\bigskip

\noindent
Let
\begin{eqnarray*}
\;\tilde{\;\omega _{g}} &=&\omega _{g}+\partial
\overline{\partial} \phi
>0\ \\
\widetilde{h} &=&he^{-\phi }
\end{eqnarray*}

\noindent
Let $\widetilde{h}_{m}$ be the induced Hermitian metric
of $\widetilde{h}$ on $L^{m}$, $\{\widetilde{S}
_{0}^{m},\widetilde{S}_{1,...,}^{m}\widetilde{S}_{d_{m}-1}^{m}\}$
be any orthonormal basis of $H^{0}(M,L^{m})$, where $d_{m}=\dim
H^{0}(M,L^{m}),$ with respect to the inner product

$$
(S_{1},S_{2})_{\widetilde{h}_{m}}=\int_{M}\widetilde{h}%
_{m}(S_{1}(x),S_{2}(x))dV_{\widetilde{g}}\;.
$$
By Theorem 2.1, we have
\begin{eqnarray*}
\sum\limits_{i=0}^{d_{m}-1}||\widetilde{S}_{i}^{m}(x)||_{\widetilde{h}
_{m}}^{2}
&=&\widetilde{a}_{0}(x)m^{n}+\widetilde{a}_{1}(x)m^{n-1}+
\widetilde{a}_{2}(x)m^{n-2}+... \\
&=&\left( \sum\limits_{i=0}^{d_{m}-1}||\widetilde{S}
_{i}^{m}(x)||_{h_{m}}^{2}\right) e^{-m\phi }.
\end{eqnarray*}
Thus
$$
\phi =\frac{1}{m}\log \left(
\sum\limits_{i=0}^{d_{m}-1}||\widetilde{S}
_{i}^{m}(x)||_{\widetilde{h}_{m}}^{2}\right) -\frac{1}{m}\log
\left(
\widetilde{a}_{0}(x)m^{n}+\widetilde{a}_{1}(x)m^{n-1}+\widetilde{a}
_{2}(x)m^{n-2}+...\right)
$$
As $m\rightarrow +\infty$, we obtain
\begin{eqnarray*}
&&\frac{1}{m}\log \left( \widetilde{a}_{0}(x)m^{n}+\widetilde{a}
_{1}(x)m^{n-1}+\widetilde{a}_{2}(x)m^{n-2}+...\right) \\
&=&\frac{1}{m}\log
m^{n}(1+\widetilde{a}_{1}(x)m^{-1}+\widetilde{a}
_{2}(x)m^{-2}+...) \\
&=&\frac{n}{m}\log m+\frac{1}{m}\log (1+O(\frac{1}{m}))\rightarrow
0
\end{eqnarray*}

\noindent Thus we have the following corollary of the
Tian-Yau-Zelditch expansion.

\begin{corollary}
$$
\left\| \phi -\frac{1}{m}\log \left(
\sum\limits_{i=0}^{d_{m}-1}||\widetilde{
S}_{i}^{m}(x)||_{h_{m}}^{2}\right) \right\| _{C^{k}}\rightarrow
0,\;as\;m\rightarrow +\infty.
$$
\end{corollary}
In other words, any plurisubharmonic function can be approximated
by the logarithms of holomorphic sections of $L^{m}$.

\bigskip
\bigskip
\centerline{\bf \S 3. Proof of Theorem 1}
\bigskip
\setcounter{theor}{3}
\noindent Let $M$\ be the blow-up of
$CP^{2}$\ at two points and $\pi $\ its natural projection.
Without loss of generality, we may assume the two points are
$p_{1}=[0,1,0]$ and $p_{2}=[0,0,1].$ Then $M$ is a subvariety of $
CP^{2}\times CP^{1}\times CP^{1}$ defined by the equations
$$
Z_{0}X_{1}=Z_{1}X_{0},\;\;Z_{0}Y_{2}=Z_{2}Y_{0}
$$
where $Z_{i},$\ $X_{j},Y_{k}$ are respectively the homogeneous
coordinates on $CP^{2}$, $CP^{1}$ and $CP^{1}$.

\bigskip
\noindent Let $G$\ be the automorphism group acting on
$CP^{2}\times CP^{1}\times CP^{1} $ generated by $\theta _{j}$\
and permutations $\tau $\ $(0 \leq i \leq 2)$
$$
\theta _{j}:[Z_{0},Z_{j},Z_{2}]\times \lbrack X_{0},X_{1}]\times
\lbrack Y_{0},Y_{2}]\rightarrow \lbrack Z_{0},Z_{j}e^{i\theta
},Z_{2}]\times \lbrack X_{0},X_{1}]\times \lbrack Y_{0},Y_{2}]
$$
for $\theta \in [0, 2\pi)$, and
$$ \tau :[Z_{0},Z_{1},Z_{2}]\times \lbrack
X_{0},X_{1}]\times \lbrack Y_{0},Y_{2}]\rightarrow \lbrack
Z_{0},Z_{2},Z_{1}]\times \lbrack Y_{0},Y_{2}]\times \lbrack
X_{0},X_{1}]\;.
$$

\bigskip
\noindent Let $\pi _{0},\pi _{1},\pi _{2}$\ be the projection from
$CP^{2}\times CP^{1}\times CP^{1}$ onto $CP^{2}$, $CP^{1}$ and
$CP^{1}$. Define $\omega$ by
\begin{eqnarray*}
 \omega&=&\pi _{0}^{\ast }\omega _{0}+\pi _{1}^{\ast }\omega _{1}+\pi
_{2}^{\ast }\omega _{2} \\
&=&\partial \overline{\partial }\log
(|Z_{0}|^{2}+|Z_{1}|^{2}+|Z_{2}|^{2})+\partial \overline{\partial
}\log (|X_{0}|^{2}+|X_{1}|^{2})+\\
&& \partial \overline{\partial }\log (|Y_{0}|^{2}+|Y_{2}|^{2})
\end{eqnarray*}
where $\omega _{0}$, $\omega _{1}$, $\omega _{2}$ are the
Fubini-Study metrics $ CP^{2}$, $CP^{1}$ and $CP^{1}$. By explicit
calculation, it can be shown that $\omega |_{M}$ is in the first
Chern class of $M$(see \cite{A}).

\bigskip
\noindent Consider the divisor
$$\{[0,Z_{1},Z_{2}]\times CP^{1}\times
CP^{1}\}+\{CP^{2}\times [1,0]\times CP^{1}\}+\{CP^{2}\times
CP^{1}\times  [1,0]\}$$ which defines a line bundle $(L,h)$ on
$CP^{2}\times CP^{1}\times CP^{1}$, where
$$h=\frac{1}{
(|Z_{0}|^{2}+|Z_{1}|^{2}+|Z_{2}|^{2})(|X_{0}|^{2}+|X_{1}|^{2})(|Y_{0}|^{2}+|Y_{2}|^{2})
},$$ then $(L,h)|_{M}\rightarrow M$ defines the anticanonical line
bundle on $M$ whose curvature form  $-\partial \overline{\partial
}\log h$  gives the first Chern class on $M.$

\noindent Since $M\backslash \{\pi ^{-1}\{p_{1}\}\cup \pi
^{-1}\{p_{2}\}\}$ is isomorphic to $CP^{2}\backslash
\{p_{1},p_{2}\}$, if we choose the inhomogeneous coordinates
$(z_{1},z_{2})=[1,z_{1},z_{2}]$\ on $CP^{2}$, the K\"ahler metric
$$
\omega _{g_{0}}=\partial \overline{\partial }\log
(1+|z_{1}|^{2}+|z_{2}|^{2})+\partial \overline{\partial }\log
(1+|z_{1}|^{2})+\partial \overline{\partial }\log (1+|z_{2}|^{2})
$$
can be extended to a K\"ahler metric $g_{0}$\ on $M$\ which
belongs to $ c_{1}(M)$. If we take different inhomogeneous
coordinates $ (w_{0},w_{1})=[w_{0},w_{1},1]$, the corresponding
K\"ahler metric is
$$
\omega _{g_{1}}=\partial \overline{\partial }\log
(1+|w_{0}|^{2}+|w_{1}|^{2})+\partial \overline{\partial }\log
(1+|w_{0}|^{2})+\partial \overline{\partial }\log
(|w_{0}|^{2}+|w_{1}|^{2})
$$
and we have

\begin{eqnarray*}
\det g_{0} &=&\frac{1}{(1+|z_{1}|^{2}+|z_{2}|^{2})^{3}}+\frac{1}{
(1+|z_{1}|^{2}+|z_{2}|^{2})^{2}(1+|z_{1}|^{2})} \\
&&+\frac{1}{(1+|z_{1}|^{2}+|z_{2}|^{2})^{2}(1+|z_{2}|^{2})}+\frac{1}{
(1+|z_{1}|^{2})^{2}(1+|z_{2}|^{2})^{2}} \\
\det g_{1} &=&\frac{1}{(1+|w_{0}|^{2}+|w_{1}|^{2})^{3}}+\frac{1}{
(1+|w_{0}|^{2}+|w_{1}|^{2})^{2}(|w_{0}|^{2}+|w_{1}|^{2})} \\
&=&\frac{1}{(1+|w_{0}|^{2}+|w_{1}|^{2})^{2}(1+|w_{0}|^{2})}+\frac{|w_{0}|^{2}
}{(1+|w_{0}|^{2})^{2}(|w_{0}|^{2}+|w_{1}|^{2})^{2}}.
\end{eqnarray*}

\bigskip
\noindent Consider the line bundle $(L^{N},h_{N})\rightarrow
CP^{2}\times CP^{1}\times CP^{1}$. Then
$$\dim H^{0}(CP^{2}\times CP^{1}\times
CP^{1},O(L^{N}))=\frac{(N+1)^{3}(N+2)}{2}$$ and
$\{Z_{0}^{i_{0}}Z_{1}^{i_{1}}Z_{2}^{i_{2}}X_{0}^{j_{0}}X_{1}^{j_{1}}
Y_{0}^{k_{0}}Y_{2}^{k_{2}}\}_{i_{0}+i_{1}+i_{2}=j_{0}+j_{1}=k_{0}+k_{2}=N,}
$ is an orthogonal basis for $H^{0}(CP^{2}\times CP^{1}\times
CP^{1},O(L^{N}))$.

\noindent Let $M_{1}$ be the hypersurface of $CP^{2}\times
CP^{1}\times CP^{1}$\ defined by the equations $$
Z_{0}X_{1}=Z_{1}X_{0} $$ and $M_{2}$\ be the hypersurface of
$CP^{2}\times CP^{1}\times CP^{1}$\ defined by the equations $$
Z_{0}Y_{2}=Z_{2}Y_{0}. $$ Then $M=M_{1}\cap M_{2}.$

\noindent In view of the short exact sequences
\begin{eqnarray*}
0 &\rightarrow &O(L^{N}-[M_{1}])\rightarrow O(L^{N})\rightarrow
O(L^{N}|_{M_{1}})\rightarrow 0 \\
0 &\rightarrow &O(L^{N}|_{M_{1}}-[M])\rightarrow
O(L^{N}|_{M_{1}})\rightarrow O(L^{N}|_{M})\rightarrow 0
\end{eqnarray*}
we can choose $N$\ sufficiently large so that
\begin{eqnarray*}
H^{1}(CP^{2}\times CP^{1}\times CP^{1},O(L^{N}-[M_{1}]))
=H^{1}(M_{1},O(L^{N}|_{M_{1}}-[M]))=0.
\end{eqnarray*}
Then $ H^{0}(CP^{2}\times CP^{1}\times CP^{1},O(L^{N}))\rightarrow
H^{0}(M_{1},O(L^{N}|_{M_{1}}))\rightarrow 0 $
$$ H^{0}(M_{1},O(L^{N}|_{M_{1}}))\rightarrow
H^{0}(M,O(L^{N}|_{M}))\rightarrow 0 $$ Hence
$\{Z_{0}^{i_{0}}Z_{1}^{i_{1}}Z_{2}^{i_{2}}X_{0}^{j_{0}}X_{1}^{j_{1}}Y_{0}^{k_{0}}Y_{2}^{k_{1}}|_{M}\}_{i_{0}+i_{1}+i_{2}=j_{0}+j_{1}=k_{0}+k_{2}=N}
$ contains an orthogonal basis for $H^{0}(M,O(L^{N}|_{M}))$ and
\begin{eqnarray*}
||Z_{0}^{i_{0}}Z_{1}^{i_{1}}Z_{2}^{i_{2}}X_{0}^{j_{0}}X_{1}^{j_{1}}Y_{0}^{k_{0}}Y_{2}^{k_{1}}||_{h_{N}}^{2}
=\frac{\left|
Z_{0}^{i_{0}}Z_{1}^{i_{1}}Z_{2}^{i_{2}}Z_{0}^{j_{0}}Z_{1}^{j_{1}}Z_{0}^{k_{0}}Z_{2}^{k_{2}}\right|
^{2} }{
(|Z_{0}|^{2}+|Z_{1}|^{2}+|Z_{2}|^{2})(Z_{0}|^{2}+|Z_{1}|^{2})(Z_{0}|^{2}+|Z_{2}|^{2}))^{N}}
\end{eqnarray*}
on $CP^{2}\backslash \{p_{1},p_{2}\}$. By Corollary 2.1, for any
$\varphi$ in $P_{G}(M,\omega_{g})$, we have on $CP^{2}\backslash
\{p_{1},p_{2}\}$,

\begin{eqnarray*}
&&\varphi ([Z_{0},Z_{1},Z_{2}])\\
&=&\lim_{N\rightarrow \infty }\frac{1}{N}\log
\frac{\sum\limits_{_{i_{0}+i_{1}+i_{2}=j_{0}+j_{1}=k_{0}+k_{2}=N}}^{{}}|a_{(
\varphi
)i_{0}i_{1}i_{2}j_{0}j_{1}k_{0}k_{2}}^{(N)}Z_{0}^{i_{0}+j_{0}+k_{0}}Z_{1}^{i_{1}+j_{1}}Z_{2}^{i_{2}+k_{2}}|^{2}}{((|Z_{0}|^{2}+|Z_{1}|^{2}+|Z_{2}|^{2})(|Z_{0}|^{2}+|Z_{1}|^{2})(|Z_{0}|^{2}+|Z_{2}|^{2}))^{N}}
\end{eqnarray*}
for some coefficients  $a_{(\varphi
)i_{0}i_{1}i_{2}j_{0}j_{1}k_{0}k_{2}}^{(N)}$= $a_{(\varphi
)i_{0}i_{2}i_{1}k_{0}k_{2}j_{0}j_{1}}^{(N)}$ which is set in the
view of the group action by $G$.

\bigskip

\begin{lemma}$\frac{1}{n}\log \frac{\sum
\limits_{_{i_{0}+i_{1}+i_{2}=j_{0}+j_{1}=k_{0}+k_{2}=n}}^{{}}|Z_{0}^{i_{0}+j_{0}+k_{0}}Z_{1}^{i_{1}+j_{1}}Z_{2}^{i_{2}+k_{2}}|^{2}%
}{
((|Z_{0}|^{2}+|Z_{1}|^{2}+|Z_{2}|^{2})(|Z_{0}|^{2}+|Z_{1}|^{2})(|Z_{0}|^{2}+|Z_{2}|^{2}))^{n}
}\leq Const$.
\end{lemma}

\begin{proof}

\noindent On the patch $U_{0}=\{Z_{0}\neq 0\}$, let
$z_{1}=\frac{Z_{1}}{Z_{0}}$ and $z_{2}=\frac{Z_{2}}{Z_{0}}$,

\begin{eqnarray*}
&&\frac{1}{n}\log \frac{\sum
\limits_{_{i_{0}+i_{1}+i_{2}=j_{0}+j_{1}=k_{0}+k_{2}=n}}|Z_{0}^{i_{0}+j_{0}+k_{0}}Z_{1}^{i_{1}+j_{1}}Z_{2}^{i_{2}+k_{2}}|^{2}
}{
((|Z_{0}|^{2}+|Z_{1}|^{2}+|Z_{2}|^{2})(|Z_{0}|^{2}+|Z_{1}|^{2})(|Z_{0}|^{2}+|Z_{2}|^{2}))^{n}} \\
&\leq& \frac{1}{n}\log \left(
\sum\limits_{_{i_{0}+i_{1}+i_{2}=j_{0}+j_{1}=k_{0}+k_{2}=n}}\frac{
|z_{1}^{i_{1}+j_{1}}z_{2}^{i_{2}+k_{2}}|^{2}}{
(1+|z_{1}|^{2}+|z_{2}|^{2})^{n}(1+|z_{1}|^{2})^{n}(1+|z_{2}|^{2})^{n}}\right)
\\
&\leq&\frac{1}{n}\log \left(
\sum\limits_{_{i_{0}+i_{1}+i_{2}=j_{0}+j_{1}=k_{0}+k_{2}=n}}^{{}}\frac{
|z_{1}^{i_{1}+j_{1}}z_{2}^{i_{2}+k_{2}}|^{2}}{
1+|z_{1}^{i_{1}+j_{1}}z_{2}^{i_{2}+k_{2}}|^{2}}\right)\\
&\leq& \frac{1}{n}\log \left(
\sum\limits_{_{i_{0}+i_{1}+i_{2}=j_{0}+j_{1}=k_{0}+k_{2}=n}}^{{}}1\right) \\
&=&\frac{1}{n} \log \frac{(n+1)^{3}(n+2)}{2}
\end{eqnarray*}

\noindent On the patch $U_{2}=\{Z_{2}\neq 0\}$, let
$w_{0}=\frac{Z_{0}}{Z_{2}}$ and $w_{1}=\frac{ Z_{1}}{Z_{2}},$

\begin{eqnarray*}
&&\frac{1}{n}\log \frac{\sum
\limits_{_{i_{0}+i_{1}+i_{2}=j_{0}+j_{1}=k_{0}+k_{2}=N}}^{{}}|Z_{0}^{i_{0}+j_{0}+k_{0}}Z_{1}^{i_{1}+j_{1}}Z_{2}^{i_{2}+k_{2}}|^{2}
}{
((|Z_{0}|^{2}+|Z_{1}|^{2}+|Z_{2}|^{2})(|Z_{0}|^{2}+|Z_{1}|^{2})(|Z_{0}|^{2}+|Z_{2}|^{2}))^{n}
} \\
&\leq&\frac{1}{n}\log \left(
\sum\limits_{_{i_{0}+i_{1}+i_{2}=j_{0}+j_{1}=k_{0}+k_{2}=n}}^{{}}\frac{
|w_{0}^{i_{0}+j_{0}+k_{0}}w_{1}^{i_{1}+j_{1}}|^{2}}{
(1+|w_{0}|^{2}+|w_{1}|^{2})^{n}(1+|w_{0}|^{2})^{n}(|w_{0}|^{2}+|w_{1}|^{2})^{n}
}\right) \\
&\leq&\frac{1}{n}\log \left(
\sum\limits_{_{i_{0}+i_{1}+i_{2}=j_{0}+j_{1}=k_{0}+k_{2}=n}}^{{}}\frac{
|w_{0}^{i_{0}+j_{0}+k_{0}}w_{1}^{i_{1}+j_{1}}|^{2}}{
|w_{0}^{i_{0}+j_{0}+k_{0}}w_{1}^{i_{1}+j_{1}}|^{2}}\right) <
\frac{1}{n} \log \left(
\sum\limits_{_{i_{0}+i_{1}+i_{2}=j_{0}+j_{1}=k_{0}+k_{2}=n}}^{{}}1\right) \\
&=&\frac{1}{n}\log \frac{(n+1)^{3}(n+2)}{2}
\end{eqnarray*}

\noindent This inequality holds for the patch $U_{1}=\{Z_{1}\neq
0\}$, and so the lemma is proved.
\end{proof}

\begin{lemma} There exists $\varepsilon>0$ such that for any $\varphi\in P_{G}(M,\omega_{g})$ and N,
there exist $n>N$, $i_{0},i_{1},i_{2},j_{0},j_{1},k_{0},k_{2}$
with $i_{0}+i_{1}+i_{2}=j_{0}+j_{1}=k_{0}+k_{2}=n$, and
$(a_{(\varphi
)i_{0}i_{1}i_{2}j_{0}j_{1}k_{0}k_{2}}^{(n)})^{\frac{1}{n}}>\varepsilon$.
\end {lemma}

\begin{proof}

\noindent Otherwise, for any $\varepsilon >0,$ there exists $
\varphi $ and $N,$\ such that for any $\ n>N$ and any $
i_{0},i_{1},i_{2},j_{0},j_{1},k_{0},k_{2}$ satisfying
$i_{0}$+$i_{1}$+$i_{2}$=$j_{0}$+$j_{1}$=$k_{0}$+$k_{2}$=$n$, we
have $(a_{(\varphi
)i_{0}i_{1}i_{2}j_{0}j_{1}k_{0}k_{2}}^{(n)})^{\frac{1}{n}}<\varepsilon$.
By choosing $n$ large enough and with the lemma above, we have

\begin{eqnarray*}
&&\varphi ([Z_{0},Z_{1},Z_{2}])\\
&\leq &\frac{1}{n}\log \frac{\max |a_{(\varphi
)i_{0}i_{1}i_{2}j_{0}j_{1}k_{0}k_{2}}^{(n)}|^{2}\sum
\limits_{_{i_{0}+i_{1}+i_{2}=j_{0}+j_{1}=k_{0}+k_{2}=n}}^{{}}|Z_{0}^{i_{0}+j_{0}+k_{0}}Z_{1}^{i_{1}+j_{1}}Z_{2}^{i_{2}+k_{2}}|^{2}
}{
((|Z_{0}|^{2}+|Z_{1}|^{2}+|Z_{2}|^{2})(|Z_{0}|^{2}+|Z_{1}|^{2})(|Z_{0}|^{2}+|Z_{2}|^{2}))^{n}}+\varepsilon
\\
&\leq&\frac{1}{n}\log \frac{\sum
\limits_{_{i_{0}+i_{1}+i_{2}=j_{0}+j_{1}=k_{0}+k_{2}=n}}^{{}}|Z_{0}^{i_{0}+j_{0}+k_{0}}Z_{1}^{i_{1}+j_{1}}Z_{2}^{i_{2}+k_{2}}|^{2}}{
((|Z_{0}|^{2}+|Z_{1}|^{2}+|Z_{2}|^{2})(|Z_{0}|^{2}+|Z_{1}|^{2})(|Z_{0}|^{2}+|Z_{2}|^{2}))^{n}}
+2\log\varepsilon+\varepsilon
\\
&\leq &\log \varepsilon +const
\end{eqnarray*}

\noindent Since $\varepsilon $\ could be arbitrarily small, the
above inequality would imply that $\varphi \rightarrow -\infty $\
uniformly, which contradicts the fact that $\sup_{M}\varphi =0.$
\end{proof}

\bigskip
\noindent \textbf{Proof of the theorem:}

\begin{eqnarray*}
&&\varphi ([Z_{0},Z_{1},Z_{2}])\\
&=&\lim_{N\rightarrow \infty }\frac{1}{N}\log
\frac{\sum\limits_{_{i_{0}+i_{1}+i_{2}=j_{0}+j_{1}=k_{0}+k_{2}=N}}^{{}}|a_{(
\varphi
)i_{0}i_{1}i_{2}j_{0}j_{1}k_{0}k_{2}}^{(N)}Z_{0}^{i_{0}+j_{0}+k_{0}}Z_{1}^{i_{1}+j_{1}}Z_{2}^{i_{2}+k_{2}}|^{2}
}{
((|Z_{0}|^{2}+|Z_{1}|^{2}+|Z_{2}|^{2})(|Z_{0}|^{2}+|Z_{1}|^{2})(|Z_{0}|^{2}+|Z_{2}|^{2}))^{N}
} \\
&\geq &\frac{1}{N}\log \frac{
|Z_{0}^{i_{0}+j_{0}+k_{0}}Z_{1}^{i_{1}+j_{1}}Z_{2}^{i_{2}+k_{2}}|^{2}+|Z_{0}^{i_{0}+j_{0}+k_{0}}Z_{1}^{i_{2}+k_{2}}Z_{2}^{i_{1}+j_{1}}|^{2}
}{
((|Z_{0}|^{2}+|Z_{1}|^{2}+|Z_{2}|^{2})(|Z_{0}|^{2}+|Z_{1}|^{2})(|Z_{0}|^{2}+|Z_{2}|^{2}))^{N}
}+C_{1} \\
&\geq &\frac{1}{N}\log
\frac{|Z_{0}{}^{m}Z_{1}^{\frac{3}{2}N-\frac{m}{2}
}Z_{2}^{\frac{3}{2}N-\frac{m}{2}}|^{2}}{
((|Z_{0}|^{2}+|Z_{1}|^{2}+|Z_{2}|^{2})(|Z_{0}|^{2}+|Z_{1}|^{2})(|Z_{0}|^{2}+|Z_{2}|^{2}))^{N}
}+C_{1} \\
&\geq &\log
\frac{|Z_{0}{}^{\frac{2m}{N}}Z_{1}^{3-\frac{m}{N}}Z_{2}^{3-\frac{m}{N}}|^{2}}{
(|Z_{0}|^{2}+|Z_{1}|^{2}+|Z_{2}|^{2})(|Z_{0}|^{2}+|Z_{1}|^{2})(|Z_{0}|^{2}+|Z_{2}|^{2})
}+C_{1}
\end{eqnarray*}
where $i_0+j_0+k_0=m,i_1+j_1+i_2+k_2=3N-m$.

\bigskip
\noindent On the patch $U_{0}=\{Z_{0}\neq 0\}$,
\begin{eqnarray*}
&& \int_{U_{0}\cap \{0< \left| z_{1}\right| ,\left| z_{2}\right| <
1\}}e^{-\alpha \varphi }\omega _{g_{0}}^{2} \\
&\leq& C_{1}\int_{0< \left| z_{1}\right| ,\left| z_{2}\right| <
1}e^{-\alpha \log
\frac{|Z_{0}|^{\frac{2m}{N}}|Z_{1}|^{3-\frac{m}{N}}|Z_{2}|^{3-\frac{m}{N}}}{
(|Z_{0}|^{2}+|Z_{1}|^{2}+|Z_{2}|^{2})(|Z_{0}|^{2}+|Z_{1}|^{2})(|Z_{0}|^{2}+|Z_{2}|^{2})
}}\omega _{g_{0}}^{2} \\
&=& C_{1}\int_{0< \left| z_{1}\right| ,\left| z_{2}\right| <
1}\frac{ (|Z_{0}|^{2}+|Z_{1}|^{2}+|Z_{2}|^{2})^{\alpha
}(|Z_{0}|^{2}+|Z_{1}|^{2})^{\alpha
}(|Z_{0}|^{2}+|Z_{2}|^{2})^{\alpha }}{ |Z_{0}|^{\frac{2\alpha
m}{N}}|Z_{1}|^{3\alpha -\frac{\alpha m}{N}
}|Z_{2}|^{3\alpha -\frac{\alpha m}{N}}}\omega _{g_{0}}^{2} \\
&\leq& C_{2}\int_{0< \left| z_{1}\right| ,\left| z_{2}\right| < 1}
\frac{(1+|z_{1}|^{2}+|z_{2}|^{2})^{\alpha }(1+|z_{1}|^{2})^{\alpha
}(|1+|z_{2}|^{2})^{\alpha }}{|z_{1}|^{3\alpha -\frac{\alpha m}{N}
}|z_{2}|^{3\alpha -\frac{\alpha m}{N}}}dz_{1}\wedge
d\overline{z}_{1}\wedge
dz_{2}\wedge d\overline{z}_{2} \\
&\leq& C_{3}\int_{0< \left| z_{1}\right| ,\left| z_{2}\right| < 1}
\frac{1}{|z_{1}|^{3\alpha -\frac{\alpha m}{N}}|z_{2}|^{3\alpha
-\frac{\alpha
m}{N}}}dz_{1}\wedge d\overline{z}_{1}\wedge dz_{2}\wedge d\overline{z}_{2} \\
& \leq& C_{3}\int_{0< \left| z_{1}\right| ,\left| z_{2}\right| <
1} \frac{1}{|z_{1}|^{3\alpha }|z_{2}|^{3\alpha }}dz_{1}\wedge
d\overline{z} _{1}\wedge dz_{2}\wedge d\overline{z}_{2}
\end{eqnarray*}

\bigskip
\noindent On the patch $U_{2}=\{Z_{2}\neq 0\}$,
\begin{eqnarray*}
&& \int_{U_{1}\cap 0< \left| w_{0}\right| ,\left| w_{1}\right|\leq
1}e^{-\alpha \varphi }\omega _{g_{1}}^{2} \\
&\leq& C_{4}\int_{0< \left| w_{0}\right| ,\left| w_{1}\right| \leq
1}e^{-\alpha \log
\frac{|Z_{0}|^{\frac{2m}{N}}|Z_{1}|^{3-\frac{m}{N}
}|Z_{2}|^{3-\frac{m}{N}}}{
(|Z_{0}|^{2}+|Z_{1}|^{2}+|Z_{2}|^{2})(|Z_{0}|^{2}+|Z_{1}|^{2})(|Z_{0}|^{2}+|Z_{2}|^{2})
}}\omega _{g_{1}}^{2} \\
&=&C_{4}\int_{0< \left| w_{0}\right| ,\left| w_{1}\right| \leq
1}\frac{ (1+|w_{0}|^{2}+|w_{1}|^{2})^{\alpha
}(1+|w_{0}|^{2})^{\alpha
}(|w_{0}|^{2}+|w_{1}|^{2})^{\alpha }}{|w_{0}|^{\frac{2\alpha m}{N}}|w_{1}|^{3\alpha -\frac{\alpha m}{N}}}\omega _{g_{1}}^{2} \\
&\leq& C_{5}\int_{0< \left| w_{0}\right| ,\left| w_{1}\right| \leq
1} \frac{(1+|w_{0}|^{2}+|w_{1}|^{2})^{\alpha
}(1+|w_{0}|^{2})^{\alpha }(|w_{0}|^{2}+|w_{1}|^{2})^{\alpha
}}{|w_{0}|^{\frac{2\alpha m}{N} }|w_{1}|^{3\alpha -\frac{\alpha
m}{N}}(|w_{0}|^{2}+|w_{1}|^{2})}dw_{0}\wedge
d\overline{w}_{0}\wedge dw_{1}\wedge d\overline{w}_{1} \\
&\leq& C_{6}\int_{0< \left| w_{0}\right| ,\left| w_{1}\right| \leq
1} \frac{1}{|w_{0}|^{\frac{2\alpha m}{N}}|w_{1}|^{3\alpha
-\frac{\alpha m}{N} }(|w_{0}|^{2}+|w_{1}|^{2})^{1-\alpha
}}dw_{0}\wedge d\overline{w}_{0}\wedge
dw_{1}\wedge d\overline{w}_{1} \\
&\leq& C_{6}\int_{t=0}^{1}\int_{s=0}^{1}\frac{1}{s^{\frac{\alpha
m}{N}}t^{
\frac{3}{2}\alpha -\frac{\alpha m}{2N}}(s+t)^{1-\alpha }}dsdt \\
&\leq& C_{6}\int_{s=0}^{1}\frac{1}{s^{\frac{\alpha
m}{N}}t^{\frac{3}{2} \alpha -\frac{\alpha m}{2N}}s^{(1-\alpha
)p}t^{(1-\alpha )q}}dsdt
\end{eqnarray*}

\noindent where $(p+q=1)$.

\bigskip
\noindent Case 1: \ If $\ 1\leq \frac{m}{N}\leq 3$, then
\begin{eqnarray*}
\frac{\alpha m}{N}+(1-\alpha )p &<&1\Leftarrow \alpha
<\frac{1-p}{3-p}<
\frac{1-p}{\frac{m}{N}-p} \\
3\alpha -1 &<&1 \\
\frac{3}{2}\alpha -\frac{\alpha m}{2N}+(1-\alpha )q &<&1\Leftarrow
\alpha <1<\frac{1-q}{\frac{3}{2}-\frac{m}{2N}-q}
\end{eqnarray*}

\noindent Case 2: \ If $\ 0<\frac{m}{N}<1$, then
\begin{eqnarray*}
\frac{\alpha m}{N}+(1-\alpha )p &<&1\Leftarrow \alpha <1 \\
3\alpha -1 &<&1 \\
\frac{3}{2}\alpha -\frac{\alpha m}{2N}+(1-\alpha )q &<&1\Leftarrow
\alpha < \frac{1-q}{\frac{3}{2}-q}
\end{eqnarray*}

\noindent So we could choose any $\alpha <\frac{1}{3}$, which
implies that $\alpha _{G}(M,\omega )\geq \frac{1}{3}$.

\bigskip
\bigskip
\noindent Conversely, we choose
\begin{eqnarray*}
\varphi _{\varepsilon } &=&\log( \frac{ |Z_{0}|^{6}}{
(|Z_{0}|^{2}+|Z_{1}|^{2}+|Z_{2}|^{2})(|Z_{0}|^{2}+|Z_{1}|^{2})(|Z_{0}|^{2}+|Z_{2}|^{2})
}+\varepsilon)\\
&&-\log (1+\varepsilon )\\
 &\in& P_{G}(M,\omega )
\end{eqnarray*}

\noindent Then we have $\sup_{M}\varphi _{\varepsilon }=0$\ and
$\varphi _{\varepsilon }=\log \frac{\varepsilon }{1+\varepsilon }$
on the exceptional divisors. Furthermore, we have
$$
\lim_{\varepsilon \rightarrow 0}\int_{M}e^{-\alpha \varphi
_{\varepsilon }}\omega ^{2}=\infty ,\;for\;any\;\alpha
>\frac{1}{3}.
$$
Hence we have shown $\alpha _{G}(M,\omega )=\frac{1}{3}$.

\bigskip
\bigskip
\centerline{\bf \S 4. Proof of Theorem 2}
\bigskip

\newcounter{theor1}
\setcounter{theor1}{4}
\newtheorem{theorem4}{Theorem}[theor1]
\newtheorem{lemma4}{Lemma}[theor1]
\newtheorem{corollary4}{Corollary}[theor1]
\newtheorem{claim4}{Claim}[theor1]
\newtheorem{proposition4}{Proposition}[theor1]

\noindent In this section, we will prove the generalized
Moser-Trudinger inequality on any K\"ahler manifold $M$ of
dimension $n$ whose $\alpha(M)$ is greater than $\frac{n}{n+1}$.
The following theorem is due to Tian and Zhu \cite{TZ}.
\bigskip
\begin{theorem4}
Let $(M,\omega) $ be a K\"ahler-Einstein manifold with Ric$\left(
\omega \right) =\omega ,$ then there exist constants $\delta
=\delta (n)$ and $C=C\left( n,\lambda _{2}(\omega )-1\right) \geq
0$ such that for any $\phi \in P(M,\omega )$ which satisfies $\phi
\perp \Lambda _{1}$, we have
$$
\ F_{\omega }(\phi )\geq J_{\omega }(\phi )^{\delta }-C,
$$
which is the same as
$$
\frac{1}{V}\int_{M}e^{-\phi }\omega ^{n}\leq Ce^{J_{\omega }(\phi )-\frac{1}{%
V}\int_{M}\phi \omega ^{n}-J_{\omega }(\phi )^{\delta }}.
$$
\end{theorem4}
\bigskip
\noindent This implies in particular the Moser-Trudinger
inequality on $S^2$, which reads
$$\frac{1}{4\pi}\int_{S^2}e^{-\phi }\omega \leq e^{\frac{1}{8\pi}\int_{S^2}|\nabla\phi|^2\omega
  -\frac{1}{4\pi}\int_{S^2}\phi}
$$
\bigskip
\noindent
For any $\phi\in P(M,\omega)$, put
$\omega'=\omega_{\phi}=\omega+
\partial\overline{\partial}\phi$ and $Ric(\omega)=\omega+\partial\overline{\partial}h_{\omega}$. Consider
the Monge-Amp\`ere equation
$$(\omega'+\partial\overline{\partial}\psi)^n=e^{h_{\omega}-t\psi}\omega'^n$$
\noindent We will use the continuity method backwards and let
${\phi_{t}}$ be a smooth family which solve the above equation.

$\bigskip $

\noindent The following lemmas are well-known \cite{T4}, but we
add the proofs for the sake of completeness.
\begin{lemma4}
$Ric(\omega _{t})\geq t\omega _{t},$ and we have equality if and
only if \ $t=1.$
\end{lemma4}

\begin{proof}

\begin{eqnarray*}
Ric(\omega _{t}) &=&-\partial \overline{\partial }\log \omega
_{t}^{n} =-\partial \overline{\partial }\log \frac{\omega
_{t}^{n}}{\omega ^{n}} +Ric(\omega )  =-\partial
\overline{\partial }\left( h_{\omega }-t\phi _{t}\right)
+\omega +\partial \overline{\partial }h_{\omega } \\
&=&\omega +t\phi _{t} =t\omega _{t}+(1-t)\omega  \geq t\omega
_{t}.
\end{eqnarray*}
\end{proof}

\begin{lemma4}For any $\phi \in P(M,\omega )$, if the Green's function of $\ \omega'=\omega +\partial
\overline{\partial }\phi $ is bounded from below, we have:
$$
-\inf_{M}\phi \leq \frac{1}{V}\int_{M}(-\phi )\omega'^{n}+C.
$$
\end{lemma4}

\begin{proof}
\noindent Since $\omega +\partial \overline{\partial }\phi
=\omega'$ and $\omega' -\partial \overline{\partial }\phi >0$, we
have $ \Delta _{\omega ' }\phi \leq n $.
\begin{eqnarray*}
-\phi &=&\frac{1}{V}\int_{M}(-\phi )\omega'^{n}+\frac{1}{V}
\int_{M}\Delta _{\omega'}\phi (y)G_{\omega ' }(x,y)\omega'^{n} \\
&\leq &\frac{1}{V}\int_{M}(-\phi )\omega'^{n}+\frac{1}{V}
\int_{M}n(G_{\omega ' }(x,y)-\inf G_{\omega ' }(x,y))\omega'^{n} \\
&\leq &\frac{1}{V}\int_{M}(-\phi )\omega '^{n}+C.
\end{eqnarray*}
\end{proof}

\noindent Let \ $\left( M,\omega \right) $ be a K\"ahler-Einstein
manifold with Ric$ \left( \omega \right) =\omega$ and let
$P(M,\omega ,K)$ = $\{\phi \in P(M,\omega )\mid G_{\omega
+\partial \overline{\partial }\phi }(x,y)\geq -K)$. Then we have:

\begin{proposition4}
Let \ $\left( M,\omega \right) $ be a K\"ahler-Einstein manifold
with $Ric \left( \omega \right) =\omega .$ If $\alpha
(M)>\frac{n}{n+1},$\ then there exist constants $\delta
(n,\alpha,K)$\ and $C(n,\alpha,\lambda_{2}(\omega)-1,K)$ such that
for any $ \phi \in P(M,\omega ,K),$ we have
$$
\ F_{\omega }(\phi )\geq \delta J_{\omega }(\phi )-C.
$$
\end{proposition4}

\begin{proof}
Let $\omega ' =\omega +\partial \overline{\partial } \phi ,$ where
$\phi \in P(M,\omega ,K).$

\begin{eqnarray*}
\frac{1}{V}\int_{M}e^{-\alpha \phi }\omega ^{n} &=&\frac{1}{V}
\int_{M}e^{-(\alpha _{1}+\alpha _{2}+\varepsilon )\phi }\omega ^{n} \\
&\leq &\frac{1}{V}\int_{M}e^{-(\alpha _{1}+\alpha _{2})\phi
}\omega ^{n}e^{-\varepsilon \inf_{M}\phi },
\end{eqnarray*}
take \ $p=\frac{1}{\alpha _{1}},q=\frac{1}{1-\alpha _{1}}$, we
have

\begin{eqnarray*}
\frac{1}{V}\int_{M}e^{-(\alpha _{1}+\alpha _{2})\phi }\omega ^{n}
&\leq & \frac{1}{V}(\int_{M}e^{-\alpha _{1}p\phi }\omega
^{n})^{1/p}(\int_{M}e^{-\alpha _{2}q\phi }\omega ^{n})^{1/q} \\
&=&\frac{1}{V}(\int_{M}e^{-\phi }\omega ^{n})^{\alpha
_{1}}(\int_{M}e^{-
\frac{\alpha _{2}}{1-\alpha _{1}}\phi }\omega ^{n})^{1-\alpha _{1}} \\
&\leq &Ce^{\alpha _{1}J_{\omega }(\phi )-\frac{\alpha
_{1}}{V}\int_{M}\phi \omega ^{n}}(\int_{M}e^{-\frac{\alpha
_{2}}{1-\alpha _{1}}\phi }\omega ^{n})^{1-\alpha _{1}}
\end{eqnarray*}
by Lemma 4.2,
\begin{eqnarray*}
e^{-\varepsilon \inf_{M}\phi } &\leq &e^{\frac{\varepsilon }{V}
\int_{M}(-\phi )\omega '^{n}+C} \\
&=&e^{\varepsilon I_{\omega }(\phi )-\frac{\varepsilon
}{V}\int_{M}\phi
\omega ^{n}+C} \\
&\leq &e^{\varepsilon (n+1)J_{\omega }(\phi )-\frac{\varepsilon
}{V} \int_{M}\phi \omega ^{n}+C}.
\end{eqnarray*}
By Holder inequality,
\begin{eqnarray*}
\frac{1}{V}\int_{M}e^{-\phi }\omega ^{n} &\leq &(\frac{1}{V}%
\int_{M}e^{-\alpha \phi }\omega ^{n})^{\frac{1}{\alpha }} \\
&\leq &Ce^{\frac{\alpha _{1}+(n+1)\varepsilon }{\alpha }J_{\omega }(\phi )-%
\frac{\alpha _{1}+\varepsilon }{\alpha V}\int_{M}\phi \omega
^{n}}(\int_{M}e^{-\frac{\alpha _{2}}{1-\alpha _{1}}\phi }\omega ^{n})^{\frac{%
1-\alpha _{1}}{\alpha }} \\
&=&Ce^{\frac{\alpha _{1}+(n+1)\varepsilon }{\alpha }J_{\omega }(\phi )-\frac{%
1}{V}\int_{M}\phi \omega ^{n}+\frac{\alpha _{2}}{V}\int_{M}(\phi
-\sup \phi )\omega ^{n}}(\int_{M}e^{-\frac{\alpha _{2}}{1-\alpha
_{1}}(\phi -\sup \phi
)}\omega ^{n})^{\frac{1-\alpha _{1}}{\alpha }} \\
&\leq &Ce^{\frac{\alpha _{1}+(n+1)\varepsilon }{\alpha }J_{\omega }(\phi )-%
\frac{1}{V}\int_{M}\phi \omega ^{n}}(\int_{M}e^{-\frac{\alpha
_{2}}{1-\alpha _{1}}(\phi -\sup \phi )}\omega
^{n})^{\frac{1-\alpha _{1}}{\alpha }}
\end{eqnarray*}
We need to determine $\alpha _{1,}$ $\alpha _{2,}$ $\varepsilon ,$
which satisfy the following conditions
\begin{eqnarray*}
\alpha &=&\alpha _{1}+\alpha _{2}+\varepsilon >1 \\
\alpha &>&\alpha _{1}+(n+1)\varepsilon \\
1 &>&\alpha _{1}
\end{eqnarray*}
So we will choose
\begin{eqnarray*}
\alpha _{2} &=&n\varepsilon +\varepsilon ' \\
\alpha _{1} &=&1-\alpha _{2}-\varepsilon +\varepsilon
"=1-(n+1)\varepsilon -\varepsilon ' +\varepsilon "
\end{eqnarray*}
where \ $\varepsilon$ , $\varepsilon '$ , $\varepsilon "<<1,$ and
$ \varepsilon ' =o(\varepsilon )$, $\varepsilon "=o(\varepsilon '
)$.

\noindent Since $\alpha (M)>\frac{n}{n+1},$  then we can choose
$\varepsilon$, $\varepsilon '$ , $\varepsilon "$ small enough,
then we have
$$
\frac{\alpha _{2}}{1-\alpha _{1}}=\frac{n\varepsilon +\varepsilon' }{%
(n+1)\varepsilon +\varepsilon ' -\varepsilon "}<\alpha (M)
$$
and
$$
\ \int_{M}e^{-\frac{\alpha _{2}}{1-\alpha _{1}}(\phi -\sup \phi
)}\omega ^{n}<Const.
$$
Combined with the inequalities above, we have
$$ \frac{1}{V}\int_{M}e^{-\phi }\omega ^{n}\leq Ce^{(1-\delta
)J_{\omega }(\phi )-\frac{1}{V}\int_{M}\phi \omega ^{n}}
$$
Which proves the lemma.
\end{proof}

\noindent \textbf{Proof of Theorem 2}

\noindent We assume $\omega $ is the K\"ahler-Einstein metric of \
$M.$ For any $\phi \in $\ \ $P(M,\omega ),$ put $\omega ' =\omega
+\partial \overline{\partial }\phi .$ Consider $ (\omega '
+\partial \overline{\partial }\psi )=e^{h_{\omega ' }+t\psi }. $
By solving the Monge-Amp\`ere equation backwards, we get the
solutions $\phi _{t}, $ and $\phi _{1}=-\phi .$

\noindent
For $t>\frac{1}{2},$ let \ $\omega _{t}=\omega ' +\partial \overline{%
\partial }\phi _{t}=$\ $\omega +\partial \overline{\partial }(\phi _{t}-\phi
_{1}), $ by Lemma 4.1,
$$
Ric(\omega _{t})\geq \frac{1}{2}\omega _{t}.
$$
which shows that the Green function of \bigskip $\omega _{t}$ is
uniformly bounded from below. Thus by proposition 4.1 and the
calculation in \cite{TZ} we have
\begin{eqnarray*}
F_{\omega }(\phi _{t}-\phi _{1}) &\geq &\delta J_{\omega }(\phi
_{t}-\phi
_{1})-C \\
&\geq &C_{1}osc_{M}(\phi _{t}-\phi _{1})-C_{2}
\end{eqnarray*}
and consequently,
\begin{eqnarray*}
n(1-t)J_{\omega }(\phi ) &=&n(1-t)J_{\omega ' }(\phi _{1}) \\
&\geq &(1-t)(I_{\omega ' }(\phi _{1})-J_{\omega ' }(\phi _{1})) \\
&\geq &F_{\omega ' }(\phi _{t})-F_{\omega ' }(\phi _{1}) \\
&=&F_{\omega }(\phi _{t}-\phi _{1}) \\
&\geq &C_{1}osc_{M}(\phi _{t}-\phi _{1})-C_{2}
\end{eqnarray*}

\begin{eqnarray*}
F_{\omega }(\phi ) &=&-F_{\omega ' }(-\phi ) \\
&=&\int_{0}^{1}(I_{\omega ' }(\phi _{t})-J_{\omega ' }(\phi
_{t}))dt \\
&\geq &(1-t)(I_{\omega ' }(\phi _{t})-J_{\omega ' }(\phi _{t})) \\
&\geq &\frac{1-t}{n}J_{\omega ' }(\phi _{t}) \\
&\geq &\frac{1-t}{n}J_{\omega ' }(\phi
_{1})-2(1-t)(C_{1}osc_{M}(\phi
_{t}-\phi _{1})-C_{2}) \\
&\geq &\frac{1-t}{n}J_{\omega }(\phi )-2(1-t)^{2}nC_{1}J_{\omega
}(\phi )-C_{3}
\end{eqnarray*}

\noindent The theorem follows by choosing
$(1-t)<\frac{1}{2n^{2}C_{1}}.$

\end{document}